\input amstex
\magnification=1200

\bigskip
\centerline{\bf Continuity of the Volume of Simplices in Classical Geometry}

\bigskip
\centerline{\bf Feng Luo}

\bigskip
\noindent
\centerline{\bf Abstract} 

\it It is proved that the volume of spherical or hyperbolic simplices, when considered as a function
of the dihedral angles,  can be extended
continuously to degenerated simplices.\rm

\bigskip
\noindent
\S 1. {\bf Introduction}
\medskip
\noindent
1.1. It is well known that the area of a spherical or a hyperbolic 
triangle can be expressed as an affine function of the inner angles
by the Gauss-Bonnet formula. In particular, the area considered as a function of the
inner angles can be 
extended continuously to degenerated spherical or hyperbolic
triangles. 
The purpose of the paper is to show that the continuous extension 
property holds in any dimension.
Namely, if a sequence of  spherical (or hyperbolic) n-simplices has the property that their corresponding
dihedral angles at codimension-2 faces converge, then the volumes of the simplices converge.
Note that if we consider the area as a function of the three edge 
lengths of a triangle, then there does not exist any
continuous extension of the area to all degenerated triangles. For 
instance, a degenerated spherical 
triangle of edge lengths $0, \pi, \pi$ is represented geometrically as 
the intersection of two great
circles at the north and the south poles. However, its area depends on 
the intersection angle of these two geodesics 
and cannot be defined in terms of the lengths. This 2-dimensional 
simple phenomenon still holds in high dimension for
both spherical and hyperbolic simplices.
\medskip

To state our result, let us introduce some notations. Given an n-simplex with vertices $v_1, ..., v_{n+1}$,
the i-th codimension-1 face is defined to be the (n-1)-simplex with vertices $v_1,$$ ..., $$v_{i-1},$$ v_{i+1},$$
 ..., v_{n+1}$.
The dihedral angle between the i-th and j-th codimension-1 faces is denoted by $a_{ij}$.
As a convention, we define $a_{ii} =\pi$ and call the symmetric matrix $[a_{ij}]_{(n+1) \times (n+1)}$
the \it angle matrix \rm of the simplex. It is well known that the angle matrix 
$[a_{ij}]_{(n+1) \times (n+1)}$ determines the simplex up to isometry in spherical and hyperbolic geometry.

Let $\bold R^{ m \times m}$ be the space of all real $m \times m$ matrics. Our main result is the following.

\bigskip
\noindent
{\bf Theorem 1.1.} \it Let  $X_n(1)$ and $ X_n(-1) \subset \bold R^{(n+1) \times (n+1)}$ 
be the 
spaces of angle matrices of all  n-dimensional spherical  
and hyperoblic simplices respectively.
The volume 
function $V: X_{n}(k) \to \bold R$
can be extended continuously to the closure of $X_n(k)$ in $\bold 
R^{(n+1) \times (n+1)}$ for $k=1, -1$. \rm
\bigskip

Note that both spaces $X_n(1)$ and $X_n(-1)$ are fairly explicitly known. Topologically, both of them
are homeomorphic to the Euclidean space of dimension $n(n+1)/2$.
We do not know if Theorem 1.1 can be generalized to convex polytopes of the same combinatorial type
in the 3-sphere or the hyperbolic 3-space.

The proof of the theorem for spherical simplieces is quite simple. It 
is an easy consequence of the continuity
of the function which sends a semi-positive definite symmetric matrix 
to its square root.  
The case of the hyperbolic simplices is
more subtle. It uses the continuity of the square roots of 
semi-positive definite symmetric matrices and the following
property of hyperbolic simplices.  We use $B_R(x)$ to denote the ball of radius $R$ centered at $x$.

\bigskip
\noindent
{\bf Theorem 1.2.} \it For any $\epsilon >0$ and any $r >0$, there is 
$R=R(\epsilon, r, n)$ so that for any hyperbolic n-simplex $\sigma$,
if $x \in \sigma$ is a point whose distance to each totally geodesic codimenison-1
hypersurface
containing a codimension-1 face
is at most $r$, then the volume of $\sigma -B_R(x)$ is at most 
$\epsilon$.
 \rm

\bigskip

Recall that the \it center \rm and the \it radius \rm of a simplex are defined to be the center and the radius of its inscribed 
ball.
The raduis of a hyperbolic n-simplex is well known to be  uniformaly bounded from above.
Applying Theorem 1.2 to the center of the n-simplex, we conclude that for any $\epsilon >0$, there is $R=R(\epsilon)$ so that
the volume of $\sigma -B_R(c)$  is less than $\epsilon$ for any  hyperbolic n-simplex $\sigma$ with
center $c$.

\bigskip
Recent work of [MY] produces an explicit formula expressing volume of 
spherical
and hyperbolic tetrahedra in terms of
the dihedral angles using dilogarithmic function. It is not clear if Theorem 1.1 in dimension 3 follows from
their explicit formula.

\bigskip
\noindent
1.2. Using the work of Aomoto [Ao] and Vinberg [Vi], one may express the volume of a simplex in terms of an
integral related to the Gaussian distribution (see (2.3) and (2.7)). To state Theorem 1.1 in terms of matrices, let us
introduce some notations.
For an $n \times n$ matrix $A$, we use $ad(A)$ to denote the adjacency matrix of $A$. The transpose of $A$
is denoted  by $A^t$. The ij-th entry of $A$ is denoted by $A_{ij}$.
We use $A >0$ to denote the condition that all entries in $A$ are positive. Evidently, if a matrix
$A$ is positive definite, or $ad(A) >0$, then the following function $F$ is well defined,

$$ F(A) = \sqrt{ | det(ad(A))|} \int_{R^n_{\geq 0}} e^{ -x^t ad(A) x} dx \tag 1.1$$
where $x \in \bold R^n$ is a column vector, $\bold R_{\geq 0}$ is the set of all non-negative numbers and $dx$ is the
Euclidean volume form.  Theorem 1.1 is equivalent to the following,

\bigskip
\noindent
{\bf Theorem 1.3.} \it Let $\Cal X_n =\{ A \in \bold R^{ n \times n} |$
 $A^t = A$, all $A_{ii} =1$, $A$ is positive definite\} and let $\Cal Y_n =\{ A \in \bold R^{ n \times n} |$ $A^t = A$, all $A_{ii} =1$,  $ad(A) >0$, det$A <0$, and all
principal $(n-1) \times (n-1)$ submatrices of $A$ are positive definite\}. Then the function $F:\Cal X_n \cup \Cal Y_n \to \bold R$
can be extended continuously to the closure of $\Cal X_n \cup \Cal Y_n$ in $\bold R^{ n \times n}$. \rm

\bigskip
We don't know a proof of  Theorem 1.3 without using hyperbolic geometry (i.e., Theorem 1.2). 

\bigskip
\noindent
1.3. The paper is organized as follows. In \S2, we recall the basic set 
up and the Gram matrices of simplices. Also,
we prove Theorem 1.1 for spherical simplices. In \S3, we prove Theorem 1.1 
for hyperbolic simplices assuming Theorem 1.2. We prove Theorem 1.2 in 
section \S4.

\bigskip
\noindent
1.4. I would like to thank Z.-C, Han,
Daniel Ocone, Saul Schleimer for discussions. I thank Professor Nick Higham for 
directing my attention to the results on  matrices. This work is supported in part by a research councile grant from 
Rutgers University.

\bigskip

\noindent
\S2.  {\bf Preliminaries on Spherical and Hyperbolic Simplices}

\bigskip
\noindent
We recall some of the basic material related to the spherical and 
hyperbolic simplices in this section. 
In particular, we will recall the Gram matrices, the dual simplex and 
the volume formula. We also give a
proof of Theorem 1.1 for spherical simplices.
Here are the conventions and notations. Let $\bold R^m$ denote the 
m-dimensional real vector space whose elements
are column vectors.  A diagonal matrix with diagonal entries $a_{11}, ..., a_{nn}$ will be denoted by
$diag(a_{11}, ..., a_{nn})$. A diagonal matrix is \it positive \rm if all diagonal entries are positive.
The Kronecker delta is denoted by $\delta_{ij}$.
The standard inner product in $\bold R^m$ is denoted by $(x, y) = 
x^ty$. The length of a vector $x \in \bold R^m$
is denoted by $|x| = \sqrt{( x,x)}$. We use $dx = dx_1 dx_2 ...dx_m$ to 
denote the Euclidean volume element
in $\bold R^m$ and $\bold R^m _{ \geq 0} $ to denote the set $\{ (x_1, ..., x_m) \in 
\bold R^m| x_i \geq 0$ for all $i$\}.

We will make a use of the continuity of the square root of symmetric 
semi-positive definite matrix. Recall that if
$A$ is a symmetric semi-positive definite matrix, then its square root 
$\sqrt{A}$ is the symmetric semi-positive definite matrix
so that it commutes with $A$ and its square is $A$. It is well known 
that the square root matrix is unique. Furthermore,
the square root operation, considered as a self map defined on the 
space of all symmetric semi-positive definite
matrices, is continuous (theorem 6.2.37 in [HJ]). 
\bigskip
\noindent
2.1. {\bf Gram Matrices of Spherical Simplices}
\medskip
\noindent
Let $\bold R^{n+1}$ be the Euclidean space with the standard inner 
product.
The sphere $S^n$ is $\{ x \in \bold R^{n+1} | (x, x) =1\}$. A spherical 
n-simplex $\sigma^n$ has vertices $v_1, ..., v_{n+1}$
in $S^n$ so that the vectors $v_1, ..., v_{n+1}$ are linearly 
independent. 
The codimension-1 face of $\sigma^n$ opposite
$v_i$ is denoted by $\sigma^n_i$.
Let $d_{ij}$ be the spherical distance between $v_i$ and $v_j$  and 
$a_{ij}$ be the dihedral angle between
the codimension-1 faces $\sigma^n_i$ and $\sigma^n_j$ for $i \neq j$. 
Define $d_{ii}=0$ and $a_{ii} =\pi$.
Then the \it  Gram matrix \rm of $\sigma^n$ is defined to be the matrix 
$G=[ \cos (d_{ij})] =[ (v_i, v_j)]$ and the
\it angle Gram matrix \rm of the the simplex is the matrix $G^*= [ - 
\cos(a_{ij})]$. Note that both of them are
symmetric with diagonal entries being 1. The following is a well known 
fact.

\bigskip
\noindent
{\bf Lemma 2.1.} \it The Gram matrix $G$ and  the angle Gram matrix 
$G^*$ of  a simplex are related by the following formula
$$ G^* =DG^{-1}D  \tag 2.1$$ where $D$ is a positive diagonal matrix.
  
\rm

\bigskip
\noindent
{\bf Proof.} 
Let $B=[v_1, ..., v_{n+1}]$ be the
$(n+1) \times (n+1)$ matrix whose i-th column is the i-th vertex $v_i$.
Then the  Gram matrix  $G$ of
the simplex $\sigma^n$  is $B^t B$ due to the obvious formula $v_i^tv_j = (v_i, 
v_j) = \cos (d_{ij})$. 
To relate the matrix $G^*$ with $G$, we consider the \it dual simplex. 
\rm
First, find (n+1) independent vectors $w_1, ..., w_{n+1} \in \bold 
R^{n+1}$ so that $$(v_i, w_j) =\delta_{ij}. \tag 2.2$$
Define $v^*_i= w_i/|w_i|$. Then the dual simplex of $\sigma^n$ is the 
spherical simplex with
vertices $\{v_1^*, ...,v^*_{n+1}\}$. 
If we use $W=[w_1, ..., w_{n+1}]$, then (2.2) says $B^tW = Id$. In 
particular,
$  W = (B^t)^{-1}$. Thus $ W^tW = (B^tB)^{-1}$ = $G^{-1}$. However, by 
the formula $v^*_i = w_i/|w_i|$,
we see that the Gram matrix of the dual simplex is $D(W^tW)D = 
DG^{-1}D$ where $D$ is the diagonal matrix
whose ii-th entry is $|w_i|^{-1}$. On the other hand, by the definition 
of dual simplex,  the
Gram matrix of the dual is exactly the same as the  angle Gram matrix 
of $\sigma^n$. 
Namely, the spherical distance between $v_i^*$ and $v_j^*$ is $\pi - 
a_{ij}$. Thus (2.1) follows.
QED

\bigskip

The volume of the simplex $\sigma^n$ can be calculated as 
follows (see [Ao], [Vi] ). 
For the simplex $\sigma^n \subset S^n$, let the cone in $\bold R^{n+1}$ 
based at the origin over $\sigma^n$
be $K(\sigma^n) =\{ r x \in \bold R^{n+1}|   r \geq 0$ and $  x \in \sigma^n\} 
$.
Note that the linear transformation $B : \bold R^{n+1} \to \bold 
R^{n+1}$ sending the vector
$x$ to $Bx$ takes the standard basis element $e_i$ to $v_i$. In 
particular $B(\bold R^n_{\geq 0}) = K(\sigma^n)$.
Let $\mu_k =\int_0^{\infty} x^k e^{-x^2} dx$, i.e., $\mu_{2k} =  
\sqrt{\pi}(1.3....(2k-3)(2k-1)/2^{k+1}$ and $\mu_{2k+1}
= 2.4....(2k-2)(2k)/2^{k+1}$. 
Let the volume element on $S^n$ be $ds$, then the volume $V(\sigma^n)$ of the simplex 
$\sigma^n$ is given by (see [Ao], [Vi]),

$$ V(\sigma^n) =\int_{\sigma^n} ds$$
$$ = \mu_n^{-1} \int_{K(\sigma^n)} e^{-(x,x)} dx $$
$$ = \mu_n^{-1} \int_{ B(\bold R^{n+1}_{ \geq 0})} e^{-(x,x)} dx$$
$$ = \mu_n^{-1} \int_{ \bold R^{n+1}_{\geq 0}} e^{ -(By, By)} |det 
B| dy$$
$$ =\mu_n^{-1} \sqrt{|det G|} \int_{ \bold R^{n+1}_{\geq 0}} e^{ 
-y^t Gy} dy. \tag 2.3$$

\bigskip
\noindent
{\bf Lemma 2.2.} \it Let $\chi$ be the characteristic function of the 
set $\bold R^{n+1}_{\geq 0}$ in $\bold
R^{n+1}$, then the volume  $V(\sigma^n)$ of a spherical simplex $\sigma^n$ can be written as
$$ V(\sigma^n) = \mu_n^{-1} \int_{\bold R^{n+1}} e^{-(x,x)} \chi 
(\sqrt{ G^*} (x)) dx. \tag 2.4$$ \rm

\bigskip
\noindent
{\bf Proof.}  Note that since
$G =B^tB$ is positive definite, $G^{-1}$ is again symmetric and 
positive definite. Let $A =\sqrt{G^{-1}}$ be
the square root of $G^{-1}$ so that $A$ is symmetric positive definite 
and $AGA =Id$. Now make a change
of variable $y = Az$ in (2.3) where $z \in A^{-1}( \bold R^{n+1}_{\geq 
0})$. Then,
$V(\sigma^n) = \mu_n^{-1} \int_{  A^{-1}( \bold R^{n+1}_{\geq 0})} 
e^{-(z,z)} dz$.
Note that the characteristic function of $ A^{-1}( \bold R^{n+1}_{\geq 
0})$ is the same as the composition $\chi \circ A$.
Thus the volume is $$V(\sigma^n) = \mu_n^{-1} \int_{\bold R^{n+1}} e^{-(x,x)} \chi 
( A(x)) dx. $$
Finally, note that if we make a change of variable of the form $ x 
=D(y)$ where $D$ is a positive diagonal matrix,
the integral (2.3) does not change. By lemma 
1.1, we have
$A = D \sqrt{ G^*} D$ for a positive diagonal matrix $D$. Thus (2.4) holds.

\bigskip
\noindent
2.2. {\bf A Proof of Theorem 1.1 for Spherical Simplices}

\bigskip

 We give a proof of Theorem 1.1 for spherical simplices in this 
section. Let $X_n(1) $ be the space of all angle matrices   $[a_{ij}]_{(n+1) \times (n+1)}$  of
spherical n-simplices where $a_{ij}=a_{ji}$ and $a_{ii}=\pi$.
 The map sending $[a_{ij}]$ to the angle Gram matrix  $G^*= 
[-\cos(a_{ij})]$ is an embedding
of the closure of $X_n(1)$ into the space of all semi-positive 
definite, symmetric matrices whose diagonal
entries are 1.
Thus, to prove the continuity of the volume function on $X_n(1)$, by 
(2.4) it suffices to show the continuity of  the function
$W: \Cal X_n \to \bold R$ sending a matrix $A$ to
$$W(A) = \int_{R^{n+1}} e^{-(x,x)} \chi \circ  \sqrt{A}(x) dx.  \tag 2.5$$
To this end, take a sequence \{$A_m$\} in $\Cal X_n$ so that $\lim_{m \to \infty} A_m =A$ 
in $ \bold R^{(n+1) \times (n+1)}$. To establish the existence of 
$\lim_{m \to \infty} W(A_m)$, we first use the fact that the 
function sending a semi-positive
definite matrix to its square root
is continuous (theorem 6.2.37 in [HJ]).  In particular, $\sqrt{ A_m} $ 
converges to $\sqrt{A}$.

\bigskip
\noindent
{\bf Lemma 2.3.} \it Suppose $B_m$ is a convergent sequence of 
 $(n+1) \times (n+1)$ matrices so that $\lim_{m \to \infty} B_m = B$.
If  each row vector of $B$ is  none-zero, then the function $\chi 
\circ B_m$ 
converges almost everywhere to $\chi \circ B$ in $\bold R^{n+1}.$ \rm

\bigskip
Assuming this lemma, we finish the proof as follows. Since all diagonal entries of $A$ are 1,  we conclude that no row vector in $\sqrt{A}$ is zero.
Thus by the lemma,
$\chi \circ \sqrt{A_m}$ 
converges almost everywhere to $\chi \circ \sqrt{A}$ in $\bold 
R^{n+1}$. Since the integrant in $W(A)$) is bounded by the integrable function $ e^{-(y,y)}$, the dominant convergent theorem 
implies that $\lim_{m \to \infty} W(A_m)$ exists.

To prove lemma 2.3, let $R_i =\{ x \in \bold R^{n+1} | x_i =0\}$ be the 
coordinate planes.
Then $B^{-1}(R_i)$ is a proper subspace of $\bold R^{n+1}$. Indeed, if 
otherwise, say for some index $i$,
$B(\bold R^{n+1}) \subset R_i$,
then the i-th row of $B$ must be zero.  This contradicts the 
assumption. Therefore,
the Lebegue measure of $B^{-1}(R_i)$ is zero for all indices $i$. Now 
we claim for every point
$x \in \bold R^{n+1} - \cup _{i=1}^{n+1} B^{-1}(R_i)$,  the sequence 
$\chi \circ B_m(x)$ converges to $\chi \circ B(x)$.
Indeed, by the assumption, $B_m(x)$ converges to $B(x) \in \bold 
R^{n+1} - \cup_{i=1}^{n+1}R_i$. Thus 
we have $\chi( B_m(x))$ converges to $\chi(B(x))$. QED

\bigskip

The above also produced a proof of  Theorem 1.3 for the case of continuous extension of $F$ to the closure of $\Cal X_n$.

\bigskip

\noindent
2.3. {\bf Volume and Gram Matrices of Hyperbolic Simplices}

\bigskip

The (n+1)-dimensional
Minkowski space $\bold R^{n, 1}$ is $\bold R^{n+1}$ together with the
symmetric non-singular
bilinear form $<x,y> = \sum_{i=1}^n x_i y_i - x_{n+1} y_{n+1}=x^tSy$ where
$S =diag(1,1,...,1, -1)$ is an $(n+1) \times (n+1)$ diagonal matrix.  We define
the hyperboloid of two sheets to be $S(-1) =\{ x \in \bold R^{n,1}
| <x,x> = -1\}$ and the unit sphere $S(1) =\{ x \in \bold R^{n, 1}
| <x,x>=1\}$. The  space $S(-1)$ has two connected 
components.
It is well known that each of them can be
taken as a model for the n-dimensional hyperbolic space $H^n$. For 
simplicity, we take $H^n$ to be the component with positive last coordinates, i.e.,
$H^n = S(-1) \cap \{ x_{n+1} > 0\}$.  Given a vector $u \in S(1)$, let
$u^{\perp}$ be the totally geodesic codimension-1 space $\{ x \in H^n | 
< x, u> =0\}$.
The following lemma is well known (see for instance [Vi]).

\bigskip
\noindent
{\bf Lemma 2.4.} \it Suppose $u, v \in S(1) \cup S(-1)$. The following 
holds.

(1)  If $u, v \in H^n$, then $<u,v> \leq -1$ and the hyperbolic 
distance  between $u, v$ is $\cosh^{-1} (-<u, v>)$.

(2) If $u,v \in S(1)$, then $u^{\perp}$ intersects $v^{\perp}$ if and only 
if $|<u,v>| < 1$. In this case, the
dihedral angle of the intersection $u^{\perp}$, $v^{\perp}$ in the region 
$\{ x \in H^n | <x,u><x,v> \geq 0\}$ is
$ \arccos(-<u,v>)$.

(3)  If $u \in H^n$ and $v \in S(1)$, then the distance from $u$ to 
$v^{\perp}$ is $\cosh^{-1}(\sqrt{ 1 + <u,v>^2})$. \rm

\bigskip

A hyperbolic n-simplex $\sigma^n$ has vertices $v_1, ..., v_{n+1}$ in
$H^n$ so that these vectors are linearly independent in $\bold 
R^{n,1}$.
We denote the codimension-1 face of $\sigma^n$ opposite to $v_i$ by 
$\sigma^n_i$.
The hyperbolic distance between $v_i$ and $v_j$ is denoted by $d_{ij}$ 
and
the dihedral angle between $\sigma^n_i$ and $\sigma^n_j$ is denoted by 
$a_{ij}$
for $i \neq j$. As a convention, $d_{ii}=0$ and $a_{ii} =\pi$.
As in the case of spherical simplices, we define the \it Gram matrix
\rm $G$ of $\sigma^n$ to be $G=[ \cosh d_{ij}]=[-<v_i, v_j>]$ and the \it  angle Gram 
matrix
\rm of $\sigma^n$ to be $G^* = [-\cos (a_{ij})]$. Note that both of 
these matrices
are symmetric with diagonal entries $\pm 1$. 

The counterpart of lemma 2.1 holds, it is the following,

\bigskip
\noindent
{\bf Lemma 2.5.} \it Suppose $G$ and $G^*$ are the Gram matrix and the 
angle Gram matrix of a hyperbolic n-simplex, then
there is a positive diagonal matrix $D$  so that 
$$
G^* = -DG^{-1}D.$$ \rm 

\bigskip
\noindent
{\bf Proof.} By lemma 2.4,   $\cosh d_{ij} = -<v_i, v_j> = -v_i^tSv_j$.
Let $B =[v_1,..., v_{n+1}]$ be the square matrix whose i-th column
is the i-th vertex $v_i$, then by definition the Gram matrix $G$ is $- B^tSB$ where $S =diag(1,1,..,1,-1)$. To relate 
$G^*$ with $G$, we
find vectors $w_1, ..., w_{n+1}$ in $\bold R^{n, 1}$ so that
$<v_i, w_j> = \delta_{ij}$. Indeed, these vectors can be found by 
taking the matrix $W=[w_1, ..., w_{n+1}]$.
The condition $<v_i, w_j> = \delta_{ij}$ translates to the equation,  
$B^tSW=Id$, i.e., $W = S(B^t)^{-1}$.  By the construction of vertices 
$\{v_1,..., v_{n+1}\}$, the bilinear form
 $<,>$ restricted to 
the codimension-1 linear space spanned  by $ \{v_1, ..., v_{n+1} \} -\{ 
v_i  \}$
has signature $(n-1, 1)$. This implies that 
$<w_i, w_i>$ is positive.  
Define $v^*_i = w_i/ \sqrt{<w_i, w_i>}$. Then $v_i^* \in S(1)$ and 
$<v_i^*, v_j > =( \sqrt{<w_i, w_i>})^{-1} \delta_{ij}$.
The last equation shows that $v^*_i$ is the unit vector in $S(1)$ 
orthogonal to the i-th codimension-1 face $\sigma^n_i$
so that $< v_i, v_i^*>  > 0$. By lemma 2.4(2), the intersection angle 
$a_{ij}$ between  $\sigma^n_i$ and $\sigma^n_j$
is given by the equation
$-\cos a_{ij} =< v_i^*, v^*_j>$. This shows that the Gram matrix $A 
=[<v^*_i, v^*_j>]$ of the vectors $\{v_1^*, ..., v_{n+1}^*\}$ is
equal to the angle Gram matrix $G^*$. On the other hand, $v_i^* = w_i/\sqrt{ < w_i, w_i>}$. 
Thus the Gram matrix $A$ can be expressed
as $D F D$ where $D$ is a diagonal matrix with positive diagonal 
entries and $F$ is the Gram matrix $[ < w_i, w_j>]$.
By definition, $F = W^tSW$. Since $W=S(B^t)^{-1}$, we have $F= W^tSW = 
(B^tSB)^{-1} = -G^{-1}$. This establishes
$G^* = -DG^{-1} D$. QED

\bigskip
Let the volume element on $H^n$ be $ds$, let $K(\sigma^n)=\{ rx \in \bold R^{n+1} | r \geq 0, x \in \sigma^n\}$ be the cone 
based at the vertex 0 spanned by the simplex $\sigma^n$
in the vector space $\bold R^{n+1}$ and $dx = dx_1 ... dx_{n+1}$ be the 
Euclidean volume form in the Euclidean metric
in $\bold R^{n+1}$. Then the hyperbolic volume $V(\sigma^n)$ is given 
by (see [Vi], p28, note the Gram matrix used
in [Vi] is the angle Gram matrix in our case),

$$ V(\sigma^n) =\int_{\sigma^n} ds$$
$$ = \mu_n^{-1} \int_{K(\sigma^n)} e^{<x,x>} dx  $$
$$ = \mu_n^{-1} \int_{ B(\bold R^{n+1}_{ \geq 0})} e^{<x,x>} dx$$
$$ = \mu_n^{-1} \int_{ \bold R^{n+1}_{\geq 0}} e^{ <By, By>} |det 
B| dy$$
$$ =\mu_n^{-1} \sqrt{|det G|} \int_{ \bold R^{n+1}_{\geq 0}} e^{ 
y^tB^tS By} dy $$
$$ =\mu_n^{-1} \sqrt{|det G|} \int_{ \bold R^{n+1}_{\geq 0}} 
e^{-y^t Gy} dy. \tag 2.6$$  
Since the integration in (2.6) remains unchanged if we replace $G$ by $DGD$ 
for a positive diagonal matrix,
by lemma 2.5,  (2.6) is the same as
$$ V(\sigma^n) = 
\mu_n^{-1}( \sqrt{|det G^*|})^{-1} \int_{ \bold R^{n+1}_{\geq 0}} e^{ 
y^t (G^*)^{-1}y} dy $$
$$
= \mu_n^{-1} \sqrt{ |det(ad(G^*))|}     \int_{ \bold R^{n+1}_{\geq 0}} e^{ -y^t ad(G^*)y} dy \tag 2.7$$

\bigskip
To summary, we have

\bigskip 
\noindent
{\bf Lemma 2.6.} ([Vi]) \it Suppose a hyperbolic n-simplex has angle Gram matrix $G^*$. Then the volume of the simplex
is a function of $G^*$ given by (2.7). \rm

\bigskip
\noindent
2.4. {\bf Some Results from Matrix Pertubation Theory}

\bigskip
The following two results will be used frequently in the paper. See [SS], [Wi] for  proofs. The first theorem states the
continuous dependence of eigenvalues on the matrices. 

\bigskip
\noindent
{\bf Theorem 2.7}(Ostrowski) \it Let $\lambda$ be an eigenvalue of $A$ of algebraic multiplicity $m$. Then for
any matrix norm $|| .||$ and all sufficiently small $\epsilon >0$, there is $\delta >0$ so that if $|| B-A|| \leq \delta$, the
disk $\{ z \in \bold C | |z -\lambda | \leq \epsilon \}$ contains exactly $m$ eigenvalues of $B$ counted with multiplicity. \rm

\bigskip
The next theorem concerns the continuous dependence of eigenvectors on the matrices. We state the result in the
form applicable to our situation. Recall that an eigenvalue of a matrix is called \it simple \rm if it is the simple root of the
characteristic polynomial. 

\bigskip
\noindent
{\bf Theorem 2.8} (see [Wi], p67) \it Suppose $A_m$ is a sequence of  $n \times n$ matrices converging to $B$. Suppose $\lambda$ is
a simple  eigenvalue of $B$ and $\lambda_m$ is a simple eigenvalue of $A_m$ so that $\lim_{m \to \infty} \lambda_m =\lambda$.
Then there exists
a sequence of eigenvectors $v_m$ of $A_m$ associated to $\lambda_m$ so that these eigenvectors  converge to an
eigenvector of $B$ associated to $\lambda$. \rm

\bigskip
This theorem follows from the  fact that 
 if $\lambda$ is simple eigenvalue, then  the adjacency matrix $ad(B -\lambda Id)$ has rank 1 and its non-zero column vectors
are the eigenvectors of $B$ associated to $\lambda$. 

\bigskip
\noindent
\S3. {\bf A Proof of Theorem 1.1 for Hyperbolic Simplices Assuming 
Theorem 1.2}

\bigskip
\noindent
Recall that $X_n(-1)$ denotes the space of all angle matrices $[a_{ij}]$ of hyperbolic n-simplices.
The map $\cos (x)$ is an embedding of $[0, \pi]$ to $[-1, 1]$. Thus the angle Gram matrix $G^*=[-\cos (a_{ij})]$
is a map which embeds the closure of $X_n(-1)$ in $\bold R^{(n+1) \times (n+1)}$ to the space of
all symmetric matrics. The characterization of angle Gram matrix $[-\cos (a_{ij})]$ was known. 

\bigskip
\noindent
{\bf Lemma 3.1.} ([Lu], [Mi]) \it An $(n+1) \times (n+1)$ symmetric matrix $A$ with diagonal entries being one is the angle Gram matrix
of a hyperbolic n-simplex if and only if

\noindent
(3.1) all principal $n \times n$ submatrices of $A$ are positive definite, 

\noindent
(3.2) $det(A) <0$, and,

\noindent
(3.3) all entries of the adjacency matrix $ad(A)$ are positive. \rm

\bigskip
Let $\Cal Y_{n+1}$ be the space of all real matrices satisfying conditions in lemma 3.1 and define a function $F: \Cal Y_{n+1} \to
\bold R$ as in (1.1).  Note that by change the variable $x$ to $D(x)$ for a positive diagonal matrix $D$, we see that $F(A) = F(DAD)$.
Thus to establish  theorem 1.1 for hyperbolic n-simplices, it suffices to prove  that $F: \Cal Y_{n+1} \to \bold R$
can be extended continuously to the closure $\bar{ \Cal Y_{n+1}}$ in $\bold R^{(n+1) \times (n+1)}$.  This will be the goal in the
rest of the section.

\bigskip
\noindent
3.1. To prove Theorem 1.3 for $\Cal Y_{n+1}$, take a convergent sequence of matrices $A_m \in \Cal Y_{n+1}$ so that
$\lim_{m \to \infty} A_m = A_{\infty}$ where $A_{\infty} \in \bold R^{(n+1) \times (n+1)}$. We will prove 
that $\lim_{ m \to \infty} F(A_m)$ exists. Since the function $F(A) = F(DAD)$ for any positive diagonal matrix $D$,
we will modify the sequence $\{A_m\}$ by $D_m A_m A_m$ for positive diagonal matrices $D_m$
 so that $\lim_{m \to \infty}
F(D_mA_mD_m)$ converges. This will be the strategy of the proof.

By definition, all diagonal
entries of $A_{\infty}$ are 1. If $det(A_{\infty}) \neq 0$, then the signature of $A_{\infty}$ is $(n, 1)$. If $det(A_{\infty}) =0$,
we claim that $A_{\infty}$ is semi-positive definite. Indeed, by definition, all principal proper submatrices of $A_{\infty}$ are semi-positive
definite. This, together with $det(A_{\infty}) =0$, implies that $A_{\infty}$ is semi-positive definite. The proof of Theorem 1.3 
uses the following lemma to perturb $A_{\infty}$ and $A_m$ to $DA_{\infty}D$ and $D_m A_m D_m$ for some
positive diagonal matrices $D$ and $D_m$ so that $D_m A_m D_m$ converges to $DA_{\infty} D$ and all non-zero
eigenvalues of $D_m A_m D_m$ and $DA_{\infty} D$ are simple, i.e., they are the simple roots of the characteristic polynomials.

\bigskip
\noindent
{\bf Lemma 3.2.} \it Given a symmetric $n \times n$ matrix $A$ of signature $(k,0)$ or $(k, 1)$, and $\epsilon >0$, there exists a positive diagonal matrix $D$ so that $|D -Id| \leq \epsilon$ and all non-zero eigenvalues of $DAD$ are simple. \rm
\bigskip

This is a very simple consequence of the work on multiplicative inverse eigenvalue problem (see for instance [Fr]). For completeness, we provide
a simple proof of it in the appendix.

Applying this lemma, we find a positive diagonal matrix $D$ so that $DA_{\infty} D$ has only simple non-zero eigenvalues and also
a positive diagonal matrix $D_m$ within distance $1/m$ of the identity matrix so that $D_m DA_m D D_m$ has distinct eigenvalues
and $\lim_{m \to \infty} D_m DA_m D D_m = DA_{\infty} D$. Since $F(DA_mD)=F(A_m)$ for any positive diagonal matrix $D$, 
the modification of the sequence $A_m$ to $ D_m DA_m D D_m$ does not change the existence of the limit $\lim_{m \to \infty} F(A_m)$.
 By  theorems 2.7 and 2.8,  we may assume, after modifying  $A_m$ to $D_mDADD_m$, the following,

\medskip
\noindent
(3.4) all eigenvalues $\{ \lambda_i(m) | i=1,2,..., n+1\}$ of $A_m$ are pairwise distinct, i.e.,
$$ \lambda_1(m) > \lambda_2(m) > ....> \lambda_n(m) >0 > \lambda_{n+1}(m), $$
and all non-zero eigenvaules of $A_{\infty}$ are pairwise distinct.

\medskip
\noindent
(3.5) the limit $\lim_{m \to \infty} \lambda_i(m) = \lambda_i(\infty)$ exists for all $i=1,..., n+1$ where
$\lambda_i(\infty)$'s are the eigenvalues of $A_{\infty}$. Furthermore, 
either $rank(A_{\infty}) = n+1$ and 
$$\lambda_1(\infty) > \lambda_2(\infty) > ....> \lambda_n(\infty) >0 > \lambda_{n+1}(\infty),$$ 
or 
$k =rank(A_{\infty}) \leq n$ and $$\lambda_1(\infty) > \lambda_2(\infty) > ....> \lambda_k(\infty) > 
\lambda_{k+1}(\infty) =....=\lambda_{n+1}(\infty)=0.$$

\medskip
\noindent
3.2. We need the folloing canonical decomposition of matrices $A \in \Cal Y_{n+1}$.  Note that $A^2$ is symmetric and positive
definite. In particular, the symmetric positive definite matrix $B = \sqrt{ \sqrt{ A^2}}$ exists. Furthermore, the function
$B=B(A): \Cal Y_{n+1} \to \bold R^{(n+1) \times (n+1)}$ can be extended continuously to the closure $\bar{\Cal Y_{n+1}}$. Suppose
the eigenvalues of $A$ are  $ \lambda_1 \geq  \lambda_2 \geq  ....\geq  \lambda_n>0 > - \lambda_{n+1} $. Then there exists an orthonormal
matrix $U =[v_1, ... ,v_{n+1}]$ whose column vectors $v_i$ are eigenvectors of length one so that

$$ A  = U diag(\lambda_1, ..., \lambda_n, -\lambda_{n+1}) U^t. \tag 3.6$$
We can recover $B$ from (3.6) by the formula  $B = U diag(\sqrt{\lambda_1}, ..., \sqrt{\lambda_n}, \sqrt{ \lambda_{n+1}}) U^t$.
In particular, we have

$$ A = BUSU^t B \tag 3.7$$
and
$$ U^tBA^{-1}BU = S . \tag 3.8$$
Furthermore, due to $Bv_i =\sqrt{\lambda_i} v_i$,
 $$ BU=[ \sqrt{\lambda_1} v_1, ..., \sqrt{\lambda_n} v_n, \sqrt{\lambda_{n+1}} v_{n+1}].$$

Note that in general the matrix $U$ is not uniquely determined by $A$ due to the multiple eigenvalues. However, if the
eigenvalues of $A$ are pairwise distinct, then each eigenvector $v_i$ of norm 1 is determined by the associated eigenvalue
$\lambda_i$ up to sign.

\bigskip
The geometric meaning of the decomposition (3.7) is the following,

\bigskip
\noindent
{\bf Proposition 3.3.} \it Consider the hyperbolic n-simplex $\sigma=U^tB^{-1} (\bold R^{n+1}_{\geq 0}) \cap H^n$ with
codimension-i faces $\sigma_i$ for $i=1,2,..., n+1$. The point 
 $e_{n+1} = [0, ...,0,1]^t$ is in the simplex $\sigma$ and  the distance from $e_{n+1}$ to
the totally geodesic codimension-1 space $sp(\sigma_i)$ is at most $\cosh^{-1} (\sqrt{1+ \lambda_{n+1}})$ for all $i$.
 \rm

\bigskip
\noindent
{\bf Proof. } The vertices of the n-simplex $\sigma = U^tB^{-1} (\bold R^{n+1}_{\geq 0}) \cap H^n$ are
$v_i = U^tB^{-1}(e_i) / < U^tB^{-1}(e_i), U^tB^{-1}(e_i)> ^{1/2}$ where $e_i=[0,..,0,1,0,..,0]^t$ is the standard basis of $\bold R^{n+1}$.
To find the distance from $e_{n+1}$ to the codimension-1 totally geodesic space $sp(\sigma_i)$, we find the normal vector
to $sp(\sigma_i)$ as follows. Consider
the column vectors $w_1, ..., w_{n+1}$ of $W = SU^tB$. These vectors $w_i$ satisfy the conditions,

\medskip
\noindent
(3.9) $< w_i, w_i> =1$ for all $i$,

\medskip
\noindent
(3.10) $< w_i, U^tB^{-1}(e_j)> =0$ for $i \neq j$,

\medskip
\noindent
(3.11) $< w_i, U^tB^{-1}(e_i)> =1$ for all $i$.

\medskip
Indeed,  (3.9) follows from (3.7) that  $W^tSW =A$ and $A_{ii}=1$ for all $i$. Also (3.10) and (3.11) follow from the identity $W^tSU^tB^{-1} = BU S S  U^t B^{-1}
=Id$. This shows that $w_i$ is the normal vector in the de-Sitter space $S(1)$ which is perpendicular to
$sp(\sigma_i)$ so that $<w_i, U^tB^{-1}e_i> >0$.  To find the distance from $e_{n+1}$ to the codimension-1 totally
geodesic hypersurface containing a codimension-1 face,  we should calculate $< w_i, e_{n+1}>$.
Indeed, since $W^t S e_{n+1} = BUSS e_{n+1} =BUe_{n+1} = Bv_{n+1} =\sqrt{ \lambda_{n+1}} v_{n+1}$ and
 the eigenvector
$v_{n+1}$ has norm 1, we obtain $| < w_i, e_{n+1}>| \leq \sqrt{ \lambda_{n+1}}$. By lemma 2.4(3), we conclude that the
distance from $e_{n+1}$ to these codimension-1 faces are at most $\cosh^{-1}( \sqrt{ 1 + \lambda_{n+1}})$.

Finally, we need to show that $e_{n+1}$ is in the simplex $\sigma$. This is the same as showing that all entries of the eigenvector
$v_{n+1}$ have the same sign.  To this end, we need,

\bigskip
\noindent
{\bf Lemma 3.4.} \it Suppose $B$ is a symmetric $(n+1) \times (n+1)$ matrix so that all $n \times n$ principal  
submatrices in $B$ are positive definite and $det(B) \leq 0$. Then no entry in the adjacent matrix $ad(B)$ is zero. \rm

\bigskip
Assuming this lemma, we prove that all entries of the eigenvector $v_{n+1}$ have the same sign as follows. For the variable
$t \in [0, \lambda_{n+1}]$, consider the matrix $C(t) = A + t Id$.  By definition, all $n \times n$ principal submatix of $C(t)$ are
positive definite. Furthermore, $det(C(t)) \leq 0$ since the smallest eigenvalue of $C(t)$ is $t -\lambda_{n+1} \leq 0$. By the lemma,
 all entries $ad(C(t))_{ij}$ are non-zero. On the other hand, $ad(C(t))_{ij}$ is a polynomial in $t$ and is positive when $t=0$ by (3.3).
Thus all $ad(C(t))_{ij} > 0$. Now for $t=\lambda_{n+1}$, the first colume of $ad(C(\lambda_{n+1}))$ is an eigenvector of
$A$ associated to $-\lambda_{n+1}$. Since this negative eigenvalue is simple, any two associated eigenvectors are multiple of
each other. This ends the proof.

Now to prove lemma 3.4, we first note that $B ad(B) =det(B) Id$.  Also, the positive definiteness of the principal submatrice shows that
$ad(B)_{ii} >0$ for all $i$. If there is an entry $ad(B)_{ij} =0$, then $i \neq j$.
 Without loss of generality, let us assume that $ad(B)_{1 (n+1)}=0$. Let $w$ be the
first colume of $ad(B)$. The vector $w$ is not the zero vector due to  $ad(B)_{11} >0$.
 By the assumption that the principal submatrix $P$ obtained by removing the last row and column is positive definite, 
we have $w^tBw =w^tPw >0$. 
On the other hand, $Bw=det(B)[1,0,...,0]^t$ by definition
and $w^tBw = det(B) ad(B)_{11} \leq 0$ due to $det(B) \leq 0$ and $ad(B)_{11} >0$. This is a contradiction. QED

\bigskip
\noindent

\bigskip
\noindent
3.3. We now prove Theorem 1.3. Given the convergent sequence $A_m \in \Cal Y_{n+1}$ as in subsection 3.1 so that  (3.4) and (3.5) hold,
we produce a decomposition $$A_m =B_m U_m S U_m^t B_m \tag 3.12$$ as in (3.7). Let $k$ be the rank of $A_{\infty}$.
By theorem 2.8
and   (3.4) and (3.5), we may choose eigenvectors $v_1(m), ..., v_k(m)$ of unit length for $A_m$ associated to the
simple eigenvalues $\lambda_i(m)$ for 
so that
$$ \lim_{ m \to \infty} v_i(m) = v_i(\infty)$$
exists for $i=1,2,..., k$ and $v_i(\infty)$ is an eigenvector of norm 1 for  $A_{\infty}$. In particular, we see that the matrix
$$B_m U_m =[ \sqrt{\lambda_1(m)} v_1(m), ..., \sqrt{\lambda_k(m)} v_k(m), \sqrt{ \lambda_{k+1}(m)} v_{k+1}(m), ..., \sqrt{| \lambda_{n+1}(m)|} v_{n+1}(m)]$$
is converging to $[\sqrt{\lambda_1(\infty)} v_1(\infty), ..., \sqrt{\lambda_k(\infty)} v_k(\infty), 0,...,0]$ if $k \leq n$
or to $[\sqrt{\lambda_1(\infty)} v_1(\infty), $\newline$...,
\sqrt{\lambda_n(\infty)} v_n(\infty), \sqrt{ | \lambda_{n+1}(\infty)|}v_{n+1}(\infty)]$  for $k=n+1$ by (3.4) and (3.5).

Using (3.12), let us make a change of variable $x=B_mU_m(y)$ in 
$$F(A_m) = \sqrt{ | det(A_m)^{-1}|} \int_{ \bold R^{n+1}_{\geq 0}}  e^{x^t A_m^{-1} x} dx.$$
 We obtain by (3.8),

$$F(A_m) = \int_{(B_mU_m)^{-1}(\bold R^{n+1}_{\geq 0})} e^{ y^tSy} dy$$
$$ = \int_{\bold R^{n+1}} e^{<y,y>} \chi \circ (B_mU_m) (y) dy \tag 3.13$$
where $\chi$ is the characteristic function of $\bold R^{n+1}_{\geq 0}$ in $\bold R^{n+1}$.

By the contruction, $B_mU_m$ converges to a matrix in $\bold R^{(n+1) \times (n+1)}$. We claim that
the sequence of functions $\chi \circ (B_mU_m)$ converges almost everywhere in $\bold R^{n+1}$.  In fact, by lemma 2.3,
it suffices to verify that no row vector in $\lim_{m \to \infty} B_mU_m$ is zero. 
Suppose otherwise, say the i-th row is zero. Then the ii-th entry in $\lim_{m \to \infty} B_mU_mSU_m^tB_m$ is zero. But by assumption, the ii-th entry in $B_mU_mSU_m^tB_m$  is 
$(A_m)_{ii}$ which is always 1.  

To summary, we see that the integrant in (3.13) converges almost everywhere in $\bold R^{n+1}$. To prove that the limit
$\lim_{m \to \infty} F(A_m)$ exists, we will use the following well known lemma from analysis. We omit the proof.

\bigskip
\noindent
{\bf Lemma 3.5.} \it Suppose $\{ f_m\}$ is a sequence of integrable non-negative functions converging almost everywhere
to $f$ in $\bold R^n$. If for any $\epsilon >0$, there exists a measureable set $E \subset \bold R^n$ so that

(a) the restriction $f_m|_E$ converges a.e. to $f|_E$ and is dominated by an integrable function $g$ on $E$, and

(b) $\int_{E} f_m dx \leq \epsilon$ for all integer $m \geq 1$,

\noindent
then the $\lim_{m \to \infty} \int_{\bold R^n} f_m dx$ exists.  \rm

\bigskip

To apply this lemma, we will produce a decomposition of integral (3.13) as follows. For any $p>0$ and $p <1$, consider the set
$\Omega_p =\{ x \in \bold R^{n+1} | <x, x> \leq -p(x,x)\}$ where $(x,x)=x^tx$ is the Euclidean inner product.
The intersection $\Omega_p \cap H^n$ is equal to the hyperbolic ball of radius  $r=\cosh^{-1}( \sqrt{ (1+p)/2p})$ centered at $e_{n+1}$. 
Indeed, we may write $(x,x) = <x,x> +2(x, e_{n+1})^2$. Thus $<x,x> \leq -p(x,x)$ inside $H^n$ is the same as 
$|<x, e_{n+1}>| \leq \sqrt{ (1+p)/2p}$.
By lemma 2.4(1), the claim that $\Omega_p \cap H^n = B_r(e_{n+1})$ follows. Now in the region
$\Omega_p$, the integral $\int_{\Omega_p}  e^{<y,y>} \chi \circ (B_mU_m) (y) dy $ converges since the intergrant is dominated by
the integrable function $e^{-p(y,y)}$. On the other hand, the integral
$\int_{\bold R^{n+1} -\Omega_p}  e^{<y,y>} \chi \circ (B_mU_m) (y) dy $ is the same as
$\mu_n$ $vol( \sigma_m -B_r(e_{n+1}))$ where $\sigma_m = (B_mU_m)^{-1}(\bold R^{n+1}_{\geq 0})\cap H^n$ is a hyperbolic
n-simplex. By  proposition 3.3 and the existence of $\lim_{m \to \infty} \lambda_{n+1}(m)$, there is
a constant $C$ independent of $m$ so that $e_{n+1}$ is within distance $C$ to each codimension-1 totally geodesic
surface containing a codimension-1 face of the simplex $\sigma_m$. 
By Theorem 1.2 and proposition 3.3, the volume  $vol( \sigma_m -B_r(e_{n+1}))$ can be made arbitrary small for all n-simplices
$\sigma_{m}$  if the radius $r$ is large. Thus, by lemma 3.5,
we conclude that the limit $\lim_{m \to \infty} F(A_m)$ exists.  QED

\bigskip
\noindent
\S4. {\bf A Proof of Theorem 1.2}

\bigskip
We prove Theorem 1.2 in this section. Recall that $B_R(x)$ denotes
the ball of radius $R$ centered at $x$.

\bigskip
\noindent
{\bf Theorem 1.2.} \it For any $\epsilon >0$ and $r >0$, there exists a
positive number $R=R(\epsilon, r, n)$ so that for any
hyperbolic n-simplex $\sigma$, if $x \in \sigma$ is a point whose 
distance
to each totally geodesic hyperplane containing a codimension-1 face is
at most $r$, then the volume of $\sigma -B_R(x)$ is at most $\epsilon$.
\rm

\bigskip
The theorem will follow from a sequence of propositions and lemmas on
hyperbolic simplices. To begin with, we fix the notations and conventions as follow.
 The projective disk model of $H^n$ is denoted by $D^n =\{ (x_1, ..., x_n) \in \bold R^n | \sum_{i=1}^n x_i^2 < 1\}$. 
The
compact closure of $D^n$ is denoted by $\bar {D^n}$ which is the 
compactification of the hyperbolic space by adding the ideal points. The 
hyperbolic distance
in $H^n$ or $D^n$ will be denoted by $d$. If $\{v_1, ..., v_{k}\}$ is a set of 
points in $\bar{D^n}$, the  convex hull of it will be denoted by $C(v_1, 
..., v_k)$. The \it volume \rm of $C(v_1, ..., v_{n+1})$ in $D^n$, denoted by $vol(C(v_1, ..., v_k))$,
 is the hyperbolic volume of $C(v_1, ..., v_{n+1}) \cap D^n$. If $v_1, ..., v_{n+1}$ are pairwise distinct, 
we call
$\sigma =C(v_1, ..., v_{n+1})$ a \it generalized n-simplex \rm in 
$ \bar{D^n}$. Its i-th codimension-1 face, denoted by $\sigma^i$ is
$C(v_1, ..., v_{i-1}, v_{i+1}, ..., v_{n+1})$.   A 
generalized n-simplex is said to be \it non-degenerated \rm if it has positive volume.
Evidently, a generalized n-simplex $C(v_1, ..., v_{n+1})$
in $\bar{D^n}$ is non-degenerated if and only if the vectors $\{v_1, ..., 
v_{n+1}\}$
are linearly independent in $\bold R^{n+1}$.
The  \it center \rm and the \it  radius \rm of a non-degenerated generalized 
n-simplex are defined
to be the center and the radius of its inscribed ball.
Given a finite set $X \in \bar{ D^n}$ so that $X$ contains at least two 
points,
the smallest complete totally geodesic submanifold containing $X$ in its closure
is deonted by $sp(X)$. 
For a measurable subset $X$ of $H^n$, or $D^n$, we use $vol(X)$ to denote the
volume of the set. If $X$ lies in a totally geodesic submanifold of 
dimension-k
$H^k$
in $H^n$, we use $vol_k(X)$ to denote the volume of $X$ in the
subspace $H^k$.

\bigskip
\noindent
4.1. We will establish the following propositions and lemmas in order to prove 
Theorem 1.2.

The first proposition generalizes a result of Ratcliffe.
\bigskip

\noindent
{\bf Proposition 4.1.} (see [Ra], theorem 11.3.2)
\it Suppose $\sigma_m = C(v_1(m), ..., v_{
n+1}(m))$ is a sequence of generalized n-simplices in $D^n$ so that
$\lim_{m \to \infty} v_i(m) = u_i$ exists in $\bar{D^n}$ for all 
$i=1,..., n+1$
and either $\{u_1, ..., u_{n+1}\}$ contains at least three points or 
$\{u_1, ..., u_{n+1}\}$ consists of two distinct points $\{p, q\}$ so that
 both sets $\{ i | u_i=p\}$ and $\{ i | u_j = q\}$ contain more than one 
point.
Then $\lim_{m \to \infty} vol(\sigma_m) = vol(C(u_1, ..., 
u_{n+1}))$. 
\rm

\bigskip
Note that Ratcliffe proved the
proposition when $C(u_1, ..., u_{n+1})$ is a non-degenerated generalized n-simplex. (In [Ra], a non-degenerated
generalized simplex in our sense is  called a generalized n-simplex.)
However,
if one exams his proof carefully in ([Ra], p527-529),  the non-degeneracy condition is never used.
Ratcliffe in fact already proved the proposition under
the assumption that $\{u_1, ..., u_{n+1}\}$ are pairwise distinct. 
Thus, it suffices to prove the proposition in the case that
the number of elements in $\{u_1, ..., u_{n+1}\}$ is at most $n$ and is
at least $2$ as specified in the proposition. This will be proved in subsection 4.3.

\bigskip
\noindent
{\bf Proposition 4.2.} \it For any $\epsilon >0$, there is a
number $\delta >0$ so that if the radius of the inscribed ball of
a hyperbolic n-simplex is less than $\delta$, the volume of the simplex
is less than $\epsilon$.

\bigskip
\noindent
{\bf Lemma 4.3.} \it For any $\delta > 0$ and $r>0$, there exists
$R = R(\delta, r, n)$ so that for any hyperbolic n-simplex $\sigma$ of
radius at least $\delta$, if $x \in \sigma$ is a point whose distance
to each  codimension-1 totally geodesic surface containing a 
codimension-1
face is at most $r$, then $d(x, c) \leq R$ where $c$ is the center of 
$\sigma$. \rm

\bigskip

Finally,  we recall the 
following
useful lemma of Thurston,

\bigskip
\noindent
{\bf Lemma 4.4}. \it Given a generalized hyperbolic n-simplex $\sigma =
C(v_1, ..., v_{n+1})$ where $n \geq 2$, let $\tau =C(v_1, ..., v_n)$ be a 
codimension-1 face of $\sigma$, then
$$ vol_n(\sigma) \leq 1/(n-1) vol_{n-1} (\tau).$$ \rm

See [Thu], chapter 6, or [Ra], p518-528, especially p528 for a proof.

\bigskip

\noindent
4.2. {\bf A Proof of Theorem 1.2}

Assuming the results above, we finish the proof of Theorem 1.2 as 
follows.
Suppose otherwise that Theorem 1.2 is not true. Then there are $\epsilon_0 >0$,
$r_0 >0$, a sequence of hyperbolic n-simplices $\sigma_m$, and
a point $x_m \in \sigma_m$ so that,

\medskip
\noindent
(4.1) The distance of $x_m$ to the totally geodesic codimension-1 
surface
containing each codimension-1 face of $\sigma_m$ is at most $r_0$,
and,

\medskip
\noindent
(4.2) $Vol( \sigma_m - B_m(x_m)) \geq \epsilon_0$.

\medskip
By proposition 4.2 and condition (4.2), we may assume that the radius
$r_m$ of
$\sigma_m$ is at least $\delta_0 >0$ for all $m$. 
By lemma 4.3 for 
$\delta_0
$ and $r_0$, we find a constant $R_0$ so that
$d(x_m ,c_m) \leq R_0$ for all $m$ where $c_m$ is the center of the
simplex $\sigma_m$. In particular, $B_{m-R_0}(c_m) \subset B_m(x_m)$.
This implies $\sigma_m - B_m(x_m) \subset \sigma_m - B_{m-R_0}(c_m)$
and $$vol( \sigma_m - B_{m-R_0}(c_m)) \geq \epsilon_0, \tag 4.3$$ for 
all $m$.

In the projective disk model $D^n$, we put the center $c_m$ to the Euclidean 
center
$0$ of $D^n$. By taking a subsequence if necessary, we may assume
that $\sigma_m = C(v_1(m), ..., v_{n+1}(m))$ where
the limit $\lim_{m \to \infty} v_i (m) = u_i$ exists in $\bar{D^n}$.

\bigskip
\noindent
{\bf Lemma 4.5.} \it Suppose $\sigma_m = C(v_1(m), ..., v_{n+1}(m))$ is a sequence of
hyperbolic n-simplices with center 0 in the projective model $D^n$ so that the
limit $\lim_{m \to \infty} v_i (m) = u_i$ exists in $\bar{D^n}$ for all $i$. If $\liminf_{m \to \infty} vol (\sigma_m) >0$,
then either $\{u_1, ..., u_{n+1}\}$ consists of at least three points, or $\{u_1, ..., u_{n+1}\} =\{p, q\}$, $p \neq q$,
so that 
both sets $\{ i |u_i =p\}$ and $\{j | u_j =q\}$ contain at least two points. \rm

\bigskip

To prove this lemma, suppose otherwise, there are two possibilities. In the first possibility,  $\{u_1, ..., u_{n+1}\}$ consists of one point $\{p\}$. 
Then for all $m$ large, the points $v_i(m)$ are close to $p$ in
the Euclidean metric in $\bar{D^n}$. If $p$ is in $D^n$, then the volume of $\sigma_m$ tends to zero which contradicts
the assumption. If $p$ is in $S^{n-1}$, then $\sigma_m$ cannot have the center to be  0 for $m$ large.
In the second possibility, we may assume that $u_2=...=u_{n+1} \neq u_1$. In this case, consider the codimension-1 face
$\sigma_m^1 = C(v_2(m), ..., v_{n+1}(m))$. This (n-1)-simplex is  close to $u_2$ for $m$ large in the Euclidean metric. Since the
face is tangent to 0, it follows that $u_2=...=u_{n+1}=0$. This implies that  the (n-1)-dimensional volume
$vol_{n-1}(\sigma_m^1)$ tends to zero. By Thurston's inequality lemma 4.4, this implies that the volume of
$\sigma_m$ tends to zero. This is again a contradiction. QED

\bigskip
Thus, by proposition 4.1 and (4.2),  the simplex $\sigma=\sigma(u_1, ..., u_{n+1})$
has positive volume. This implies that $\sigma$ is a non-degenerated
n-simplex in $D^n$ whose center is 0. 
Let $\chi_m$ and $\chi$ be the characteristic
functions of $\sigma_m$ and $\sigma$ in $\bar{D^n}$. Then by definition,
the function $\chi_m$ converges almost everywhere to $\chi$ in $D^n$.
Furthermore, by proposition 4.1, the integral $\int_{D^n} \chi_m dv$
converges to $\int_{D^n} \chi dv$ where $dv$ is the
hyperbolic volume element in $D^n$. By Fatou's lemma (see for instance [Roy], p86, problem 9), 
this implies that for any ball of radius $R$ centered at 0,
$vol( \sigma_m -B_R(0))$ converges to $vol (\sigma -B_R(0))$.
Choose $R$ so large that $vol(\sigma -B_R(0)) \leq \epsilon_0/2$.
Then for $m$ large, we have $vol( \sigma_m -B_R(0)) < \epsilon_0$. 
But this
contradicts (4.3) for $m$ large. QED

\bigskip
\noindent
4.3. {\bf A Proof of Proposition 4.1.}
\medskip

By the work of Ratcliffe [Ra], it suffices to show the proposition in two cases. In the first case,
the number of elements in the set  $\{u_1, ..., u_{n+1} \}$ is between 3 and $n$.
 In the second case, $ \{u_1, u_2,  ..., u_{n+1}\}$ consists
of two elements $\{p, q\}$, $p \neq q$,  so that
both sets $\{i |u_i =p\}$ and $\{j | u_j =q\}$ contain at least two points. The goal is to show that  
$\lim_{m -> \infty} vol(\sigma_m) =0$ in both cases.

The proposition holds for $n=2$. Indeed, in this case, $u_1, u_2, u_3$ are pairwise distinct.
Thus the result was proved by Ratcliffe. Assume from now on that $n \geq 3$.

First of all, we claim

\bigskip
\noindent
{\bf Claim. } If  $u_i=u_j$ for $i \neq j$ so that $u_i$ is in $D^n$, then $\lim_{m -> \infty} vol(\sigma_m) =0$.

\bigskip
Indeed, by lemma 4.4, we can estimate $vol(\sigma_m) \leq 1/(n-1)! vol_1( v_i(m), v_j(m))$. Now
$vol_1(v_i(m), v_j(m)) = d(v_i(m), v_j(m))$ tends to $d(u_i, u_j) =0$.  

By this claim,  we may
assume from now on that if $u_i = u_j$, $i \neq j$, then $u_i \in S^{n-1}$.

By the assumption on $\{u_1, ..., u_{n+1}\}$, we may choose four points, say $u_1, u_2, u_3, u_4$
so that $u_1=u_2$ and either $u_3=u_4 \neq u_1$, or $\{u_1, u_2, u_3, u_4\}$ consists of three points.
By lemma 4.4, we have $vol(\sigma_m) \leq 1/((n-1)...4.3) 
vol_3(C(v_1(m), v_2(m), v_3(m), v_4(m)))$.
This implies that it suffices to prove the proposition for $n=3$ which 
we will assume.

To prove the proposition, there  are two cases to be considered: case 1, all $u_i$'s are in $S^2$, and case 2, some $u_i$'s are in $D^3$.

In the first case that all  
$u_i$'s are in $S^2$,   
let $w_1(m),..., w_4(m)$ be four points in $S^2$ so that $v_1(m), 
v_3(m)$ lie
in the geodesic from $w_1(m)$ to $
w_3(m)$ and $v_2(m), v_4(m)$ lie in the geodesic from $w_2(m)$ to $ w_4(m)$. 
We choose $w_1(m)$ to be the end point in the ray from $v_3(m)$ to $v_1(m)$ and $w_2(m)$ similarly. By 
the
construction, we still have $\lim_{m \to \infty} w_i(m) = u_i$ for 
$i=1,2,3,4$. Furthermore, by the construction $C(w_1(m), ..., w_4(m))$ contains
the tetrahedron $C(v_1(m), v_2(m), v_3(m), v_4(m))$.  In particular,
 $$vol(C(v_1(m), ...,v_4(m))) \leq vol (C(w_1(m), ..., w_4(m))).$$
Now,  the volume of the ideal tetrahedra 
$C(w_1(m), ..., w_4(m))$ can be calculated from  the cross ratio of the four vertices $w_1(m), ..., w_4(m)$.
To be more precise, by [Th], 
the volume of an ideal hyperbolic tetrahedron with vertices $z_1, z_2, 
z_3, z_4 \in \bold C$ depends continuously on the cross raio $[z_1, z_2, z_3, 
z_4]=\frac{z_1-z_3}{z_1-z_4} : \frac{z_2-z_3}{z_2-z_4}$.
In particular, if the cross ratio tends to 0, 1, or $\infty$, then the 
volume tends to 0. In our case, by the assumption, we see that the cross ratio 
of $(w_1(m), w_2(m), w_3(m), w_4(m))$ tends to the cross ration of $u_1, u_2, u_3, u_4$ which is 0, 1, or $\infty$. Thus the
volume $vol_3(C(v_1(m), v_2(m), v_3(m), v_4(m)))$ tends to 0.

In the second case  that 
one of the points of $\{u_1, u_2, u_3, u_4\}$ is in $D^3$, by the above claim, we may assume that $u_1 =u_2$
is in $S^2$. Furthermore, by the claim, we may assume that $u_3 \neq u_4$ and $u_3 \in D^3$. Note that $u_4 \neq u_1$.
Let $w_1(m),..., w_4(m)$ be four points in $S^2$ constructed as in the previous paragraph.
By the construction $C(w_1(m), ..., w_4(m))$ contains
the tetrahedron $C(v_1(m), v_2(m), v_3(m), v_4(m))$. Furthermore,
we have $\lim_{m} w_1(m) =\lim_{m} w_2(m) = u_1$, and $\lim_m w_3(m) = w_3$ and $\lim_m w_4(m)=w_4$
both exist so that the cross ration of $\{u_1, u_1, w_3, w_4\}$ is 0, 1, or $\infty$.
Thus by case 1, we see that the volume of $C(w_1(m), ..., w_4(m))$ tends to zero. This in turn implies that
the volume of $C(v_1(m), ..., v_4(m))$ tends to zero. 
This finishes the proof.

\bigskip
\noindent
4.4. {\bf A Proof of Proposition 4.2}
\bigskip
\noindent
Suppose otherwise, there is $\epsilon_0 >0$,  a sequence of hyperbolic n-simplices $\sigma_m = C(v_1(m), ..., $ \newline $v_{n+1}(m))$
with center $0$ in $D^n$ so that the radius of $\sigma_m$ is at most $1/m$ and its volume $vol(\sigma_m) \geq \epsilon_0$.
By taking a subsequence if necessary, we may assume that the limit $\lim_{m \to \infty} v_i(m) = u_i$ exists in 
$\bar{ D^n}$. Then by lemma 4.5 and proposition 4.1, we conclude that $\sigma=C(u_1, ..., u_{n+1})$ is a non-degenerate
generalized n-simplex. In particular, these vectors $u_1, ..., u_{n+1}$ are linearly independent in $\bold R^{n+1}$.
On the other hand,  since the radius of $\sigma_m$ tends to zero, we see that all codimension-1 totally geodesic
surfaces $sp\{ u_1, ..., u_{i-1}, u_{i+1}, ..., u_{n+1}\}$ contain $0$. This is impossible for a non-degenerated simplex.
QED

\bigskip

\noindent
4.4. {\bf A Proof of Lemma 4.3}

\bigskip
\noindent
Suppose otherwise, there exist $\delta_0 >$, $r_0 >0$, a sequence of n-simplices $\{\sigma_m | m \in \bold Z_{\geq 1} \}$ of radius at least
$\delta_0$, and a point $x_m \in \sigma_m$ so that

\medskip
\noindent
(4.4) $x_m$ is within $r_0$ distance to each codimension-1 totally geodesic surface containing a codimension-1 face
of $\sigma_m$, and,

\medskip
\noindent
(4.5) $d(x_m, c_m) \geq m$ where $c_m$ is the center of $\sigma_m$.

\medskip
Let us put the center $c_m$ of $\sigma_m$ to be the origin 0 of $D^n$. By chosing a subsequence if necessary,
we may assume that $\sigma_m = C(v_1(m), ..., v_{n+1}(m))$ so that $\lim_{m \to \infty} v_i(m) = u_i$ exists in
$\bar{D^n}$, and $\lim_{m \to \infty} x_m = x$ also exists in $\bar{ D^n}$. Since the radius
 of $\sigma_m$ is bounded away from zero, we apply lemma 4.5 and proposition
4.1 to conclude that the simplex $\sigma=C(u_1, ..., u_{n+1})$ is non-degenerated whose center is 0. 
Since $d(x_m, 0) \geq m$,
it follows that $x$ has to be one
of the vertex, say $u_1$ of $\sigma$. Now consider the upper-half space model $U^n$ for
the hyperbolic space so that $u_1=x$ is the infinity and the totally geodesic codimension-1
surface containing $u_2, ..., u_{n+1}$ is the unit upper hemi-sphere $S^{n-1}_+=\{ (t_1, ..., t_n) \in \bold R^n
| \sum_{i=1}^n t_i^2 =1, t_n > 0\}$.    Let the center of the simplex $\sigma$ in
this model be $C$ and the point of  the shortest distance to $C$ in $S^{n-1}_+$ be $P$. We claim that
the angle $\angle PCu_1$ at $C$ is at least $\pi/2$. This follows from the Gauss-Bonnet theorem.
Let $U^2$ be the unique 2-dimensional hyperbolic
plane containing $C$ and $u_1$ so that $U^2$ is perpendicular to $S^{n-1}_+$.
Let $Q=[0,...,0,1]^t$ be the north pole in $S^{n-1}_+$. Then by the  construction, $ Q \in U^2$
and $P \in U^2$ due to the orthogonality. If $P=Q$, then the angle $\angle PCu_1$ is $\pi$. The claim follows.
If otherwise, consider the hyperbolic quadrilateral $QPC u_1$ in $U^2$. The angle of
the quadrilateral at $Q$, $P$ and $u_1$ are $\pi/2$, $\pi/2$ and $0$ respectively.  On the other hand,  since $C$
is the center of the simplex $\sigma$, the complete
geodesic from $u_1$ to $C$ intersects the hemi-sphere $S^{n-1}_+$ at some point, say R. Thus the quadrilateral
$QPC u_1$ is inside the hyperbolic triangle $\Delta u_1 Q R$ whose inner angles are $\pi/2, 0, \theta$.
In particular, the area of this triangle is less than $\pi/2$ by the Gauss-Bonnet formula. This implies that the
area of the quadrilateral $QPCu_1$ is at most $\pi/2$. By Gauss-Bonnet formula, we conclude that
the angle $\angle PCu_1$ at $C$ is at least $\pi/2$.

On the other hand, we will derive from (4.4) and (4.5) that the angle $\angle PCu_1$ is strictly less than $\pi/2$. Thus we arrive
a contradiction. To see this, let $P_m$ be the point in the totally geodesic codimension-1 surface $sp(C(v_2(m), ...,
v_{n+1}(m)))$ which is closest to the center $c_m$ of $\sigma_m$. By the construction, the limit of the
angle $\angle P_m c_m x_m$ is equal to $\angle PCu_1$. To estimate the angle $\angle P_m c_m x_m$,
consider the two-dimensional totally geodesic plane $D_m$ which contains $c_m$ and $x_m$ so that
$D_m$ is perpendicular to $sp(C(v_2(m), ..., v_{n+1}(m)))$. By the construction $P_m$ is in the plane $D_m$.
Let $R_m$ be the point in $sp(C(v_2(m), ...,
v_{n+1}(m))$ of the shortest distance to $x_m$. Then we again have $R_m$ is in $D_m$. Consider the
quadrilateral $P_m R_m x_m c_m$ in the plane $D_m$. The angles at the vertices $P_m$ and $R_m$ are $\pi/2$.
The distances $d(c_m, P_m) \geq \delta_0$, $d(x_m, R_m) \leq r_0$ and $d(c_m, x_m) \geq m$. Thus, as
$m$ becomes large, the quadrilateral  is tending to a right angled hyperbolic triangle  with one vertex at infinity
(corresponding to $R_m$ and $x_m$). There is an edge of the triangle having finite length which is at least $\delta_0$
(corresponding to the edge between $c_m$ and $P_m$).  The accue angle
at the end  point of this finite length edge is at most $\theta =\arcsin (1/ \cosh(\delta_0)) < \pi/2$ by the cosine law. Thus, 
as $m$ tends to infinity, the angle $\angle P_m c_m x_m$ tends to a number less than or equal to $\theta$.
In particular, the angle $\angle P_mc_m x_m$ is strictly less than $\pi/2$ for $m$ large. This contradicts the
previous conclusion. QED

\bigskip
\centerline{ \bf Appendix, A Proof of Lemma 3.2}

\bigskip
We give a proof of the following lemma used in the paper.
\bigskip
\noindent
{\bf Lemma 3.2.} \it Given a symmetric $n \times n$ matrix $A$ of signature $(k,0)$ or $(k, 1)$, and $\epsilon >0$, there exists a positive diagonal matrix $D$ so that $|D -Id| \leq \epsilon$ and all non-zero eigenvalues of $DAD$ are simple. \rm

\bigskip

\noindent
{\bf Proof. } For of all, it suffices to find a positive diagonal matrix $D$ so that all  non-zero eigenvalues of $DAD$ are simple. This is due
to the fact from algebraic geometry that an algebraic subvariety in $\bold R^m$ is either the whole space or has zero Lebesgue measure.
By [Fr], the set of all diagonal matrices $D$ so that $DAD$ has a non-simple non-zero eigenvalue forms an algebraic variety $X$ in $\bold R^m$.
Thus, as long as $ X \neq \bold R^m$, we can pick $D$ in $\bold R^m$ within $\epsilon$ distance to $[1,....,1]^t$ so that $ D \notin X$.

Next, we claim that it suffices to prove the lemma for  $n \times n$ matrix $A$ so that $det(A) \neq 0$. Indeed, if $k=rank(A)$, due to the
fact that $A$ is diagonal, $A$ has exactly $k$ non-zero eigenvalues counted with multiplicity and $A$ has a non-singular
principal $k \times k$ submatrix $B$ formed by $i_1, ..., i_k$-th rows and columns of $A$. For simplicity, we assume
that $B$ is formed by the first $k$ rows and columns of $A$. Then by the result for non-singular symmetric matrix, we find
a positive diagonal matrix $D_1= diag(a_1, ..., a_k)$ so that $D_1BD_1$ has $k$ distinct eigenvalues. Consider the $n \times n$
matrix $D(t)=diag(a_1, ..., a_k, t,t,..., t)$ where $t>0$. For $t$ small, by theorem 2.7, the eigenvalues of $D(t) A D(t)$ are
close to the eigenvalues
of $D_1BD_1$ and 0. Since $D_1BD_1$ has $k$ distinct non-zero eigenvalues,
this implies that $D(t)A D(t)$ has $k$ distinct non-zero eigenvalue for $t$ small.

Finally, we prove the lemma for non-singular matrices using induction on the size of the matrix. The result clearly holds for $1 \times 1$ and
$2 \times 2$ matrices. Suppose $A$ is a non-singular $n \times n$ matrix for $n \geq 3$. Let $B$ the principal submatrix of $A$ obtained
by removing the last column and the last row. Then the signature of $B$ is either $(n-2,1)$, $(n-1, 0)$ or $(n-2, 0)$. By 
the induction
hypothesis and the argument in the previous paragraph, we find a positive diagonal matrix $D_1 =diag( a_1, ..., a_{n-1})$ so that
all $(n-1)$-eigenvalues of $B$ are distinct. Let us denote the eigenvalues of $B$ by $\lambda_1 > ...> \lambda_{n-1}$.
Now consider the positive diagonal matrix $D(t) = diag(a_1, ..., a_{n-1}, t)$ for $t >0$.  For $t$ small, the eigenvalues
$\mu_1(t) \geq \mu_2(t) \geq ...\geq \mu_n(t)$ of $D(t) A D(t)$ is close to $\{ \lambda_1, ..., \lambda_{n-1}, 0\}$.
We claim that for $t$ small $D(t)A D(t)$ has $n$ distinct eigenvalues. Indeed, if $B$ is non-singular, i.e., the set
$\{ \lambda_1, ..., \lambda_{n-1}, 0\}$ consists of $n$ distinct elements, then for 
$t >0$ small, $\mu_i (t) \neq \mu_j(t)$ for $i \neq j$.
If $B$ is singular, then $B$ is semi-positive definite of rank $n-2$. Furthermore, this implies that $A$ has signature $(n-1, 1)$.
In particular, $D(t)AD(t)$ has a negative eigenvalue, i.e., $\mu_n(t) < 0$. Now for $t$ small, we conclude
that  $\mu_1(t), ..., \mu_{n-1}(t)$
are positive and are close to the set of $n-1$ distinct numbers $\{ \lambda_1, ..., \lambda_{n-2}, 0\}$ where $\lambda_i >0$. 
This implies that for $t >0$ small, the eigenvalues $\mu_1(t), ..., \mu_{n-1}(t)$ are positive and pairwsie distinct. Since
the smallest eigenvalue $\mu_n(t) <0$, we conclude that
$D(t)A D(t)$ has n distinct eigenvalues.

\bigskip
\centerline{ \bf References}

\bigskip
[Ao] Aomoto, Kazuhiko:  Analytic structure of Schläfli function. Nagoya Math J. 68 (1977), 1--16.

[Fr] Friedland, Shmuel: 
On inverse multiplicative eigenvalue problems for matrices. 
Linear Algebra and Appl. 12 (1975), no. 2, 127--137.

[HJ] Horn, Roger A.; Johnson, Charles R.: Topics in matrix analysis. Cambridge University Press, Cambridge, 1991.

[Lu] Luo, Feng: On a problem of Fenchel.  Geom. Dedicata 64 (1997), no. 3, 277--282.

[Mi] Milnor, John: The Schlaefli differential   equality. In Collected papers, vol. 1.
  Publish or Perish, Inc., Houston, TX, 1994. 

[MY] Murakami June;  Yano Masakazu: On the volume of a hyperbolic and spherical
tetrahedron,  http://www.f.waseda.jp/murakami/papers/tetrahedronrev3.pdf

[Ra]  Ratcliffe, John G.: Foundations of hyperbolic manifolds. Graduate Texts in Mathematics, 149. Springer-Verlag, New York, 1994.

[Ro]  Royden, H. L.: Real analysis. The Macmillan Co., New York; Collier-Macmillan Ltd.,
   London 1963.

[SS] Stewart, G. W.; Sun, Ji Guang:  Matrix perturbation theory. Computer Science and
   Scientific Computing. Academic Press, Inc., Boston, MA, 1990. 

[Th] Thurston, W.: Geometry and topology,  http://www.msri.org/publications/books \newline /gt3m/

[Vi] Vinberg, E. B.: Volumes of non-Euclidean polyhedra. (Russian) Uspekhi Mat. Nauk 48 (1993), no. 2(290), 
17--46; translation in Russian Math. Surveys 48 (1993), no. 2, 15--45 

[Wi]   Wilkinson, J. H.: The algebraic eigenvalue problem. Monographs on Numerical Analysis.
   Oxford Science Publications. The Clarendon Press, Oxford University Press, New York, 1988.

\bigskip
\noindent
Department of Mathematics, Rutgers University, Piscataway, NJ 08854, USA

\bigskip
\noindent
Email: fluo\@math.rutgers.edu

\end

\end

\end